\documentclass[12pt,a4paper,reqno]{amsart}
\usepackage{graphicx}
\usepackage[linktocpage=true,colorlinks,citecolor=blue,linkcolor=blue,urlcolor=blue]{hyperref}
\usepackage{amssymb}
\usepackage{cite}
\usepackage{bbm}
\usepackage{amsmath}
\usepackage{latexsym}
\usepackage{amscd}
\usepackage{tikz}
\usepackage{amsthm}
\usepackage{mathrsfs}
\usepackage{url}
\usepackage[utf8]{inputenc}
\usepackage[english]{babel}
\usepackage{amsfonts}
\usepackage{mathtools}
\usepackage{comment}
\usepackage[autostyle]{csquotes}
\usepackage[colorinlistoftodos]{todonotes}

\def\P{\mathbb P}

\def\E{\mathbb E}

\def\CA{\mathcal{A}}

\def\le {\leqslant}

\newtheorem{theorem}{Theorem}

\newtheorem*{problem}{Problem}

\theoremstyle{remark}
\newtheorem{remark}[theorem]{Remark}

\theoremstyle{definition}

\theoremstyle{remark}

\usepackage{bm,colonequals,physics}

\def\E{\mathbb E}

\title{Escaping Chaos in Random Multiplicative Functions}
\date{\today}

\author{Max Wenqiang Xu} 
\address{Courant Institute of Mathematical Sciences, 251 Mercer Street, New York 10012, USA}
\email{maxxu1729@gmail.com}
	\dedicatory{In memory of
    Jing-Run Chen on his birthday May 22nd\\
    Dedicate to Professor Loo-Keng Hua on the 115th anniversary of his birth\\Morningside Center of Mathematics Chinese Academy of Sciences on the 30th anniversary \\ Qiuzhen College of Tsinghua University on the 5th anniversary}

\begin{document}

\begin{abstract}
Let $f(n)$ be a Steinhaus random multiplicative function. Let $A\subset [1, N]$ be a finite set of integers. We show that 
\[\frac{1}{\sqrt{|A|}} \sum_{n\in A} f(n) \xrightarrow[]{d} \mathcal{CN}(0,1)\]
forces that $|A|=o(N)$. We prove that the $o(1)$ density is sharp by showing that 
for most sets $A$, and thus confirm the existence,  with density $\rho$
such that $(1-\rho)^{-1} =o((\log \log N)^{1/2})$, we have
\[
    \frac{1}{\sqrt{(1-\rho) |A|}} \sum_{n\in A} f(n) \xrightarrow{d} \mathcal{CN}(0,1).
\] 
The extra factor $\sqrt{1-\rho}$ makes a difference as long as the density $\rho>0$. 
\end{abstract}

\maketitle

\section{Introduction}
For number theorists with little background in probability theory, perhaps the first and most surprising phenomenon in the theory of random multiplicative functions is that the partial sums 
\[ \frac{1}{\sqrt{x}}\sum_{1\le n \le x} f(n)\]
do not obey a central limit theorem in the most friendly and familiar form, e.g., converging to a complex normal distribution $\mathcal{CN}(0,1)$ with mean zero and variance one, when $f(n)$ is a Steinhaus random multiplicative function. Here, $f(n)$ are not i.i.d. random variables in fact, $f(p)$ are independent and identically distributed random variables uniformly distributed on the complex unit circle, but $f(n)$ is defined to be completely multiplicative. The obvious global dependence makes it unclear if the most ``obvious'' central limit theorem would hold or not.  

It is quite nontrivial to prove that, in fact, the full partial sums do not converge to a standard complex normal distribution \cite{Harper}. More specifically, we may simply summarize the chaotic behaviors in terms of global distribution in two aspects: 

(1). Even the sizes do not match. Arguably, the most important and the foundational work in the area is Harper's theorem \cite{HarperLow}, which states that 
\begin{equation}\label{eqn: Harper}
 \E|\sum_{1\le n \le x}f(n)|\asymp \frac{\sqrt{x}}{(\log \log x)^{1/4}}.   
\end{equation}
This at least tells us that the usual $1/\sqrt{x}$ normalization is not a valid choice to get a standard complex normal distribution.

(2). It is just not Gaussian, no matter what normalization. As proved by the breakthrough work of Gorodetsky and Wong \cite{GWFinal} (and see also earlier work by S. Hardy\cite{hardy2025}), the limiting distribution may be viewed as a Gaussian twisted by a ``random variance" $V_{\infty}$, which is genuinely a random variable, and as a consequence, even with the ``correct'' normalization, the limiting distribution is not Gaussian. See also \cite{HSX} for a softer argument for the later fact. 

The deep reasons behind the ``extra'' $(\log \log x)^{1/4}$ saving factor and the appearance of $V_{\infty}$ random variance are closely related to the so-called ``Gaussian Multiplicative Chaos''. The purpose of this note, in some sense, is to escape the impact of multiplicative chaos on the global distribution. 

It has been a rich study that if we choose the support of $f$ be a subset $A\subset [1, N]$, then in what situations, we can have this naively and nicely behaved convergence, just as a sum of i.i.d. random variables
\begin{equation}\label{eqn: nice}
    \frac{1}{\sqrt{|A|}} \sum_{n\in A} f(n) \xrightarrow[]{d} \mathcal{CN}(0,1).
\end{equation} This can be viewed as the distribution is blind to the strong global dependence of $f(n)$ if you choose $A$ properly. 
Many examples have been studied, e.g., shifted primes, short intervals, and polynomial values etc. (See \cite{SoundXu, KSX, PWX}). 

A particular question we would like to address in this short note is:
\begin{problem}\label{Prob}
What is the largest size of $|A|$ with $A\subset [1, N]$, so that \eqref{eqn: nice} holds?   
\end{problem}
One may view this problem as we would like to avoid the global chaotic effects as much as we can. The trivial first example would be the set of primes. It was proved in \cite{SoundXu}, by Soundararajan and Xu, that one can take essentially the ``multiplication table set'' which gives you a set of size $N/(\log N)^{\delta/2+o(1)}$ with this special Erdős–Tenenbaum–Ford constant $\delta=0.086...$. The proof goes through exploring the Martingale structure and verifying a certain ``multiplicative energy'' type condition, and this naturally connects to the multiplication table problem. However, this is not the best record we know now. Indeed, as a simple consequence of the forthcoming work of Harper, Soundararajan, and Xu \cite{HSX}, one can see that the short interval $A=[N-M, N]$ would already have a larger size with $M \approx N/\big(\exp((\log \log N)^{1/2+o(1)})\big)$. We remark that this form of quantity is the correct threshold for several problems (see \cite{caich2024random, Xu, chang2024}), and the essential reason is the common ballot-type feature underlying them. 

In this note, we solve the problem. We first show that $A$ can not have positive density.

\begin{theorem}\label{thm: MAIN}
    Let $f$ be a Steinhaus random multiplicative function and $A=A_N\subset [1, N]$ be large. Then 
    \[     \frac{1}{\sqrt{|A|}} \sum_{n\in A} f(n) \xrightarrow[]{d} \mathcal{CN}(0,1) \implies |A| = o(N).\]
\end{theorem}
We next find the correct normalization for the positive density case.
\begin{theorem}\label{thm M2}
    Let $f$ be a Steinhaus random multiplicative function and $A=A_N\subset [1, N]$ be large with density $\rho=\rho_N$. For any $\rho$ with $(1-\rho)^{-1} =o((\log \log N)^{1/2})$, then for most sets $A$ with density $\rho$, and thus confirm the existence, we have
\begin{equation}\label{eqn: eqn general}
    \frac{1}{\sqrt{(1-\rho) |A|}} \sum_{n\in \CA} f(n) \xrightarrow{d} \mathcal{CN}(0,1).
\end{equation}
\end{theorem}

\begin{remark} We give some remarks regarding the two theorems. 
\begin{enumerate}
    \item Combining Theorem~\ref{thm: MAIN} and Theorem~\ref{thm M2}, we see that $\rho=o(1)$ is the threshold for the most ideal central limit theorem to hold, i.e., the usual normalization. 
    \item In Theorem~\ref{thm M2}, the density requirement $(1-\rho)^{-1} =o((\log \log N)^{1/2})$ is sharp, as otherwise, a larger density set would lead to the same (up to a multiplicative constant) limiting distribution to the full sum, which we know (see \cite{GWFinal, HSX, Harper}) is genuinely different from a Gaussian.
\item We also point out a related interesting feature that shows up in \cite{HSX} when studying the limiting distribution of random multiplicative functions over short intervals $[x, x+y]$. It is proven that one needs to normalize the partial sums by a factor to compensate for the multiplicative chaos effect when $y$ is sufficiently close to $x$. But the final distribution is still Gaussian! 
\end{enumerate}
    
\end{remark}
 It is also quite likely that one might first have a different guess for the Problem. Before we give the proof, we first give a heuristic argument which is treated a bit differently from our proof, but we hope that it is already convincing enough for the readers.

\subsection*{Heuristic}
 We first give a heuristic for why one may believe such  results. 
Instead of considering a deterministic subset $A$, we start with a random $A$ with density $\rho$. Let $a(n)$ be a sequence of i.i.d Bernoulli distributions with parameter $\rho$. Our random $\mathcal{A}$ is defined to be $\{1\le n \le N: a(n)=1\}$. Note that we choose $\CA$ independent of $f(n)$ and consider the product probability space of the random sum 
\[ \frac{1}{\sqrt{|\CA|}} \sum_{n\in \CA} a(n)f(n). \]
For the sum $\sum_{n\in \CA} a(n)f(n)$, it simply has mean zero and we can easily compute the variance by using the orthogonality of $f(n)$, independence between $f$ and $a$ to get that the variance is 
\[\mathbb{V}[\sum_{1\le n \le N} a(n)f(n)] = \E [\sum_{1\le n \le N} a(n)^{2}] = \rho N. \]
Now we proceed with the first trick: subtracting the mean of $a(n)$. 
\[w(n) := a(n) - \rho.\]
Then for any $n$, we have 
\[\P[w(n)= 1-\rho] = \rho, \qquad \P[w(n)= -\rho] = 1- \rho, \qquad \E[w(n)]=0. \]
And we do the decomposition
\[\sum_{1\le n \le N} a(n)f(n) =\sum_{1\le n \le N} w(n)f(n) + \rho \sum_{1\le n \le N}f(n). \]
However, by Harper's theorem, the second sum is small with high probability, say $\ll \sqrt{N}/(\log \log N)^{1/4-\epsilon}$, and it turns out this is negligible. 
So we really just need to consider the new weighted random sum $\sum_{n\le N}w(n)f(n)$. 

Note that $w(n)$ and $f(n)$ are independent of each other. To make life easy, we do a conditioning step. In fact, we can freeze either $w(n)$ or $f(n)$. Let us do it in a less common way which might make explanations easier: freeze $f(n)$ as a given number on the complex unit circle and view the conditioned random sum as a weighted sum of i.i.d $w(n)$, and this certainly converges to a conditional Gaussian distribution. However, one can compute that the (here $\widetilde{\E}$ and $\widetilde{\mathbb{V}}$ mean conditional on $f$)
\[\widetilde{\E}[ \sum_{n\le N}f(n)w(n)] = 0, \qquad \widetilde{\mathbb{V}} [\sum_{n\le N}f(n)w(n)]   = (1-\rho) \rho N. \]
 Note that both of the statistics are uniform in the choice of values of $f(n)$. Thus, we have that 
 the joint full distribution has 
 \[\frac{1}{\sqrt{(1-\rho) \rho N}} \sum_{n\in \CA} f(n) \xrightarrow{d} \mathcal{CN}(0,1).\]
Notice that a random $\CA$ has size concentrated around $\E[|\CA|] \sim \rho N$. So the normalization factor is like $\sqrt{(1-\rho)\E[|\CA|}]$. This additional multiplicative factor $\sqrt{1-\rho}$ is suggesting the new feature i.e., if $\rho=o(1)$ then it has no effect and otherwise it is a non-trivial factor.

\subsection*{Acknowledgement}
The author is grateful to Adam J. Harper for organizing the wonderful workshop ``Multiplicative Chaos in Number Theory" at Warwick University in 2025, where the author started writing this note, and to Ye Tian for hosting him at the Morningside Center of Mathematics of the Chinese Academy of Sciences, where the paper was essentially finished. The hospitality of both institutes is greatly appreciated.

The author thanks Pawel Nosal for pointing out crucial typos in the earlier version of the paper
and thanks
Adam J. Harper, K. Soundararajan, Quanyu Tang, and Liyuan Ye for their interest in this work.
The author is supported by a Simons Junior Fellowship from the Simons Foundation. 
\newpage
\section{Proof of Theorems}
We make the heuristic argument rigorous by adding some new ingredients. 
\subsection*{Proof of Theorem~\ref{thm: MAIN}} We first show that why $|A|=o(N)$ if \eqref{eqn: nice} holds for $A=A_N$. This is essentially a simple corollary of Harper's Theorem \ref{eqn: Harper}, plus some moment convergence fact that might be missed in this problem before. 

Suppose that $|A|= \rho N$ with $\rho\gg 1$ and \eqref{eqn: nice} holds. Then by using Harper's Theorem, we know that $\rho \sum_{1\le n \le N}f(n)$ is $o(\sqrt{N})$ with high probability, which is negligible compared with the normalization $\sqrt{|A|}\gg \sqrt{N}$. Thus, we derive that it must also hold that 
\begin{equation}\label{eqn: nice 2}
    \frac{1}{\sqrt{|A|}} \sum_{1\le n\le N} (1_A(n)-\rho)f(n) \xrightarrow[]{d} \mathcal{CN}(0,1).
\end{equation}
We compute the second moments 
\begin{equation}\label{eqn: small}
  \E[| \frac{1}{\sqrt{|A|}}\sum_{1\le n\le N} (1_A(n)-\rho)f(n) |^{2}] = \frac{1}{\rho N}(\rho N (1-\rho)^{2} + (1-\rho)N \rho^{2}) = 1 - \rho <1.   
\end{equation}
We are always warned that the convergence in distribution to a standard Gaussian does \textit{not} require the second moment convergence in general. (The particular counterexample is in the situation where some rare events can make the second moment very large, say larger than 1, but the distribution is still like a standard Gaussian. See, for example, the situation in the short interval case in \cite{HSX}). But our situation here is exactly the opposite. Crucially, \eqref{eqn: small} implies that the second moment is \text{strictly} smaller than one, by using the dominated convergence theorem, say, as $\rho=\rho_N$ is bounded away from zero. Thus, uniformly smaller second moments \eqref{eqn: small} contradict convergence in distribution to a Gaussian with larger second moment \eqref{eqn: nice 2}. And this completes the upper bound proof of the theorem.

\subsection*{Proof of Theorem~\ref{thm M2}}
To prove \eqref{eqn: eqn general} holds, we use a probabilistic argument, as kind of already explained in the heuristic section, together with three ingredients: Harper's Theorem, McLeish's central limit theorem (the version established in Soundararajan-Xu \cite{SoundXu}), plus a concentration treatment. 

The goal is to show that there exists a set $A\subset [N]$ such that \eqref{eqn: eqn general}. We write $1_A(n) = (1_A(n)-\rho) +\rho$, and it is enough to show the existence of such a set $A$ with
\begin{equation}\label{eqn: 3}
\frac{1}{\sqrt{(1-\rho)|A|}} \sum_{1\le n\le N}(1_A(n)-\rho) f(n) \xrightarrow[]{d} \mathcal{CN}(0,1).
\end{equation}
This is again because, by Harper's theorem, the contribution from the subtracted sum converging to zero in $L^{1}$ and thus it is negligible.
This shifting step is where we require, which is also necessary, that the density $\rho$ can not be too close to 1, i.e. $(1-\rho)^{-1} =o((\log \log N)^{1/2})$. 

We give a construction in two steps. Firstly, we show a random set satisfies the condition. 
Recall we define $a(n)$ to be a sequence of i.i.d. random variables, with  Bernoulli distribution of parameter $\rho$, i.e. $\P[a(n) =1]=\rho$ and $\P[a(n)=0] = 1-\rho$, and write $w(n)=a(n)-\rho$. The argument in the heuristic section is basically showing that for a random set $\CA : =\{n: w(n)=1-\rho\}$, the above convergence in distribution holds for a random $\CA$. To guarantee the existence of $A$, we need a concentration argument to show that actually it is not only for a random $\CA$, but for actually most deterministic sets $A$ with density $\rho$, the convergence in distribution \eqref{eqn: 3} holds. 
To do so, we choose to condition on a different probability space compared to the argument in the heuristic. Condition on $w(n)$, we show that for almost all realization of $w(n)$, 
\begin{equation}
    \frac{1}{\sqrt{(1-\rho)\rho N}} \sum_{1\le n \le N} w(n) f(n) \xrightarrow[]{d} \mathcal{CN}(0,1).
\end{equation}
To show such a result, we employ the weighted version of McLeish martingale central limit theorem, established by Soundararajan and Xu \cite{SoundXu}. The conditions that we need to verify are (where $P(m)$ denotes the largest prime factor of $m$) 
\begin{equation}\label{eqn: cross}
 \sum_{\substack{1\le n_1,n_2, m_1,m_2\le N\\ n_1n_2 =m_1m_2\\ m_1\neq n_1, m_2\neq n_2\\ P(m_1)=P(n_1)\\P(m_2)=P(n_2)}} w(n_1)w(n_2)w(m_1)w(m_2) = o\Big(  (\sum_{1\le n \le N} w(n)^{2})^{2}\Big).   
\end{equation}
and 
\begin{equation}\label{eqn: Lind}
 \sum_{\substack{1\le n_1,n_2, m_1,m_2\le N\\ n_1n_2 =m_1m_2\\ m_1\neq n_1, m_2\neq n_2\\ P(m_1)=P(n_1) =P(m_2)=P(n_2)}} w(n_1)w(n_2)w(m_1)w(m_2) = o\Big(  (\sum_{1\le n \le N} w(n)^{2})^{2}\Big).   
\end{equation}
We claim that the quantity $ \sum_{1\le n \le N} w(n)^{2} $ is concentrated around $(1-\rho)\rho N$ with high probability. To see this, we can quickly show that 
\[\mu: = \E[ \sum_{1\le n \le N} w(n)^{2}  ] \sim (1-\rho)\rho N, \quad
\E[ | \sum_{1\le n \le N} w(n)^{2} -\mu|^{2}  ] \sim \rho(1-\rho)(1-2\rho)^{2} N.\]
An application of Markov's inequality implies what we claimed. 

So, to show that the above condition holds for almost all realizations of $w(n)$, it is enough to show that  
\begin{equation}\label{eqn: concentration}
\E[|\sum_{\substack{1\le n_1,n_2, m_1,m_2\le N\\ n_1n_2 =m_1m_2\\ m_1\neq n_1, m_2\neq n_2\\ P(m_1)=P(n_1)\\P(m_2)=P(n_2)}} w(n_1)w(n_2)w(m_1)w(m_2)|^{2}] = o((1-\rho)^{4}\rho^{4}N^{4}). 
\end{equation}
And 
\begin{equation}\label{eqn: concentration 2}
\E[|\sum_{\substack{1\le n_1,n_2, m_1,m_2\le N\\ n_1n_2 =m_1m_2\\ m_1\neq n_1, m_2\neq n_2\\ P(m_1)=P(n_1) = P(m_2)=P(n_2)}} w(n_1)w(n_2)w(m_1)w(m_2)|^{2}] = o((1-\rho)^{4}\rho^{4}N^{4}). 
\end{equation}
Note that we only need to prove our theorem for the case $\rho \gg 1/\log N$ as we have seen examples, e.g., a set of primes up to $N$. So $\rho^{4}N^{4}\gg N^{4}/(\log N)^{4}$ say. Also we have $(1-\rho)^{4}\gg 1/(\log \log N)^{2+\epsilon}$. So we only need an upper bound, say $O(N^{4}/(\log N)^{4})$ for both conditions.
To verify \eqref{eqn: concentration 2}, we simply use bound $|w(n)|\le 1$ and the inner sum is simply a counting question and it is shown in \cite[equation (1.1)]{SoundXu} that this is automatically negligible. We next verify \eqref{eqn: concentration}. Given \eqref{eqn: concentration 2} holds, we may apply the triangle inequality and use divisor bounds to deduce that we only need to verify \eqref{eqn: concentration} under the condition $n_1n_2=m_1m_2$ and $\{n_1, n_2\} \neq \{m_1, m_2\}$. Expanding the square in \eqref{eqn: concentration}, and taking expectation, the only $\{(n_1, n_2, m_1, m_2; n_3, n_4, m_3, m_4)\}$ that would survive would be 

(1). $n_1, n_2, m_1, m_2$ are all distinct and we must have $(n_3, n_4, m_3, m_4)$ is a permutation of $(n_1, n_2, m_1, m_2)$.  The total amount is bounded by $N^{2+\epsilon}$ where we used the constraint $m_1m_2 =n_1n_2$ with a divisor bound.

(2). $n_1=n_2 = n, $ for some $n$  and $n \neq m_1 , n\neq m_2$ and $m_1\neq m_2$; and the $(n_3, n_4, m_3, m_4)$ is a permutation of $(n_1, n_2, m_1, m_2)$. The total counting in this case is  $\ll N^{1+2\epsilon}$ where we again used the constraint $m_1m_2 =n_1n_2$ with a divisor bound.

This is more than we need, and it proves \eqref{eqn: concentration} and thus we complete the proof of Theorem~\ref{thm: MAIN}.

\bibliographystyle{plain}
	\bibliography{main}{}
\end{document}